\documentclass[12pt,twoside]{article}
\usepackage{amsmath,amsthm}
\usepackage{amssymb,latexsym}
\usepackage{mathrsfs}
\usepackage{amscd}
\usepackage{enumerate}
\usepackage{color}
\usepackage{cite}
\usepackage{amssymb}
\usepackage[mathscr]{euscript}
\usepackage{graphicx}
\usepackage{titletoc}
\usepackage{hyperref}
\usepackage{esint}

\titlecontents{section}[0pt]{\addvspace{1pt}\filright}
              {\contentspush{\thecontentslabel\quad}}
              {}{\titlerule*[8pt]{$\cdot$}\contentspage}

\newtheorem{theorem}{Theorem}[section]
\newtheorem{corollary}[theorem]{Corollary}

\newtheorem{proposition}[theorem]{Proposition}
\newtheorem{definition}[theorem]{Definition}

\newtheorem{example}[theorem]{Example}

\theoremstyle{definition} \theoremstyle{remark}
\numberwithin{equation}{section}

\newcommand{\diam}{{\rm diam}}

\begin{document}
\title{\large {\bf RANDOM WALKS AND INDUCED DIRICHLET FORMS\\ ON COMPACT SPACES OF HOMOGENEOUS TYPE}}

\author{Shi-Lei Kong, Ka-Sing Lau and Ting-Kam Leonard Wong}

\date{}

\maketitle

\abstract We extend our study of random walks and induced Dirichlet forms on self-similar sets \cite {KLW, KL} to  compact spaces of homogeneous type $(K, \rho ,\mu)$ \cite {CW}.  A successive partition on $K$ brings a natural augmented tree structure $(X, E)$ that is Gromov hyperbolic, and the hyperbolic boundary is H\"older equivalent to $K$. We then introduce a class of transient reversible random walks on $(X, E)$ with return ratio $\lambda$. Using Silverstein's theory of Markov chains, we prove that the random walk induces an energy form on $K$ with
$$
{\mathcal E}_K [u] \asymp \iint_{K\times K \setminus \Delta}  \frac{|u(\xi) - u(\eta)|^2}{V(\xi, \eta)\rho (\xi, \eta)^\beta} d\mu(\xi) d\mu(\eta),
$$
where $V(\xi, \eta)$ is the $\mu$-volume of the ball centered at $\xi$ with radius $\rho (\xi, \eta)$, $\Delta$ is the diagonal, and $\beta$ depends on $\lambda$. In particular, for an $\alpha$-set in ${\mathbb R}^d$, the kernel of the energy form is of order $\frac{1}{|\xi-\eta|^{\alpha +\beta}}$. We also discuss conditions for this energy form to be a non-local regular Dirichlet form.

\renewcommand{\thefootnote}{}
\footnote {\indent \hspace{5mm} Na\"im kernel, reversible random walk.}
\footnote{The research is supported in part by the HKRGC grant and SFB 1283 of the German Research Council.}
\renewcommand{\thefootnote}{\arabic{footnote}}
\setcounter{footnote}{0}


\section{Introduction}
\label{sec:1}

 In discrete potential theory \cite{Dy}, it is well-known that a transient Markov chain $\{Z_n\}_{n=0}^\infty$ on a countable set $X$ converges to a limiting random variable $Z_\infty$  on the {\it Martin boundary} ${\mathcal M}$, where the union $X \cup {\mathcal M}$ is a compactification of $X$ under a topology involving the Green function. The Martin kernel $K(\cdot,\cdot)$ on $X \times {\mathcal M}$, together with the hitting distribution $\nu$ of $Z_\infty$, plays the role of the Poisson kernel in classical analysis: any function $\varphi \in L^1({\mathcal M}, \nu)$ can be extended to a harmonic function on the ``interior"~$X$ via
$$
(H \varphi) (x) = \int_{\mathcal M}\varphi (\xi) K(x, \xi) d\nu (\xi),  \qquad x \in X.
$$
Here harmonicity is defined in terms of the transition probability of the chain.

\medskip

The Markov chains we consider in this paper are always {\it reversible}: with a discrete energy form ${\mathcal E}_X$ on $X$ given by
\begin{equation} \label{eq1.3}
{\mathcal E}_X[u] = \frac 12\sum_{x, y \in X} c(x, y) |u(x) -u(y)|^2,
\end{equation}
where $c(x,y) = c(y,x) \geq 0$ and $m(x) := \sum_{y\in X} c(x,y) < \infty$, the corresponding reversible Markov chain on $X$ is the one with transition probability  $P(x, y) = c(x, y)/m(x)$. In \cite{Si}, Silverstein showed that for such a Markov chain, the induced energy form ${\mathcal E}_{\mathcal M}$ on ${\mathcal M}$, defined by
$$
 {\mathcal E}_{\mathcal M} [\varphi] =  {\mathcal E}_X [H\varphi],
$$
has the expression
\begin{equation}\label{eq1.4}
{\mathcal E}_{\mathcal M} [\varphi] = C\iint_{\mathcal M \times \mathcal M \setminus \Delta} |\varphi(\xi) - \varphi (\eta)|^2 \,\Theta (\xi, \eta)\,d\nu (\xi) d\nu (\eta),
\end{equation}
where $C > 0$ is a constant and $\Delta$ denotes the diagonal. Here $\Theta (\cdot, \cdot)$ is called the  {\it Na\"{i}m kernel} \cite{Na}. This has an analog in classical analysis \cite{FOT}: for the Brownian motion on the closed unit disk ${\mathbb D}$ started at the center and reflected on the boundary ${\mathbb T} = \partial \mathbb{D}$,  there is an induced jump process (Cauchy process) on ${\mathbb T}$ such that the corresponding energy form on ${\mathbb T}$ has the expression in \eqref{eq1.4} with $\Theta (\xi, \eta)=  \sin^{-2} (\frac{\xi-\eta}{2}) \asymp |\xi-\eta|^{-2}$. (For $f, g >0$, by $f\asymp g$, we mean  $ C^{-1} g \leq f\leq  C g$ for some $C>0$.)

\bigskip

The above theoretical setup leads to many challenging questions (e.g., see \cite{CF, Ki2, BGPW}). Among them we are interested in the following two problems:

\begin{enumerate}[(i)]
\item Study those Markov chains  that the Martin boundaries and the hitting distributions can be identified  with some familiar sets and measures respectively.

\item With the  reversible Markov chains satisfying (i), find explicit expressions  or estimates of the Na\"im kernels in \eqref{eq1.4}, and to study the associated Dirichlet forms.
\end{enumerate}

\medskip

These problems are difficult in general as it is non-trivial to estimate the Green function and the hitting probabilities. For random walks on discrete groups, the identification problem in (i) has been studied for a long time, and the reader can refer to \cite {Ka1} for an informative survey.  More recent attempts deal with self-similar sets in fractal geometry. For such a set $K$ generated by an iterated function system (IFS), there is a symbolic space (coding space) $\Sigma^*$ which gives a convenient symbolic representation of any $\xi \in K$. If the IFS satisfies the {\it open set condition} (OSC), then the representation is unique for generic points of $K$. There are a number of studies of random walks on $\Sigma^*$ such that $K$ can be identified with the Martin boundary ${\mathcal M}$ \cite{DS, JLW, Ka2, LN1, LN2, LW2}.
The most far reaching attempt is due to Kaimanovich \cite{Ka2}: he introduced a natural {\it augmented tree} for the Sierpi\'{n}ski gasket $K$ by adding horizontal edges on $\Sigma^*$  according to the intersections of the cells at each level (Sierpi\'{n}ski graph).  He showed that the Sierpi\'{n}ski graph is hyperbolic in the sense of Gromov \cite {Gr} and the hyperbolic boundary is  H\"older equivalent to $K$. In \cite{LW1, LW2}, construction of the augmented tree and identification of the boundary were extended to all self-similar sets. These works brought into play the hyperbolic structure and hyperbolic boundary which are very useful tools for studying analysis on fractals through random walks.

 \medskip

   In this respect, the current authors introduced in \cite {KLW} a class of reversible random walks $\{Z_n\}_{n=0}^\infty$ that have constant return ratio $0<\lambda <1$ on the augmented tree of a self-similar set $K$ with the OSC. We call this walk a {\it $\lambda$-natural random walk} ($\lambda$-NRW, see Section 3). By using a theorem of Ancona \cite {An}, it was shown that the Martin boundary of $\{Z_n\}_{n=0}^\infty$ can be identified with the hyperbolic boundary, and  thus the self-similar set $K$ as well. Moreover, by estimating the ever-visiting probability of $\{Z_n\}_{n=0}^\infty$ on states in $X$ and on ${\mathcal M}$ in terms of the Gromov product in the hyperbolic graph, we obtained explicit expression of the hitting distribution  as well as estimates of the Martin kernel and the Na\"im kernel. In particular, the Na\"im kernel satisfies $\Theta (\xi, \eta)  \asymp |\xi-\eta|^{-(\alpha + \beta)}$, where $\alpha$ is the Hausdorff dimension of $K$, and $\beta= \frac {\log \lambda}{\log r}$ and $r$ is the minimal contraction ratio of the IFS.

\bigskip

The main purpose of this paper is to extend the considerations in \cite{KLW, KL} to the class of {\it compact spaces of homogeneous type}, and summarize some of the results. A space of homogeneous type is a triple $(M, \rho, \mu)$ where $M$ is a set, $\rho$ is a quasi-metric and $\mu$ is a regular Borel measure with the doubling property. This class of spaces was introduced by Coifman and Weiss in \cite {CW}, and has been studied in great detail in the theories of metric measure spaces as well as $H^p$ spaces, singular integrals and analytic capacities (see \cite {AM, DH, LuS} and the references therein). This class contains the classical domains, the $\alpha$-sets \cite {JW}, and many more examples (see \cite {CW}).  In \cite {Ch}, Christ showed that a space of homogeneous type admits a partition system (dyadic cubes) which can be represented by a tree. When $M = K$ is compact, similar to the symbolic space of the self-similar sets, we can augment this tree by adding suitable horizontal edges corresponding to neighboring cells, and the  resulting graph (augmented tree) is hyperbolic. We can extend the above $\lambda$-NRW to this situation, and obtain similar identification results and estimates of the kernels. The induced energy form in this case has the estimate
$$
{\mathcal E}^{(\beta)}_K [u] \asymp \iint_{K\times K \setminus \Delta}  \frac{|u(\xi) - u(\eta)|^2}{V(\xi, \eta)\rho (\xi, \eta)^\beta} d\mu(\xi) d\mu(\eta),
$$
where $V(\xi, \eta)$ is the $\mu$-volume of the ball center at $\xi$ with radius $\rho (\xi, \eta)$, and $\beta$ depends on $\lambda$.

 \medskip

The domain of the induced energy form ${\mathcal E}^{(\beta)}_K$ is a Besov space $\Lambda^{\beta/2}_{2,2}$. By the polarization identity, ${\mathcal E}^{(\beta)}_K$ defines a non-local regular Dirichlet form if $\Lambda^{\beta/2}_{2,2}\cap C(K)$ is dense in $C(K)$ and in $\Lambda^{\beta/2}_{2,2}$.  The family $\{\Lambda^{\beta/2}_{2,2}\}_{\beta >0}$ is decreasing in  $\beta$, and in general it is trivial for large $\beta$. This defines a  critical exponent $\beta^*$. It is an important index because in the standard cases, it is the place where a local  Dirichlet form arises (which corresponds to a Laplacian) with domain $\Lambda^{\beta^*/2}_{2,\infty}$ (while $\Lambda^{\beta^*/2}_{2,2}$ contains only constant functions). For example for classical domains we have $\beta^* =2$, and for the Sierpi\'{n}ski gasket $\beta^* = \log 5/\log 2$.  We will give some discussion on this, and provide a criterion to determine the critical exponents on the p.c.f.~self-similar sets as in \cite {KL}.

\bigskip

The proofs of the results will appear separately; for the special case on the self-similar sets, the reader can refer to \cite {LW2, KLW, KL}.  The rest of the paper is organized as follows. We first recall,  in Section 2, some preliminaries about hyperbolic graphs and the boundary theory of Markov chains. In Section 3, we study the partition systems and the tree structures for compact spaces of homogeneous type, and establish the hyperbolicity of the augmented trees. In Section 4, we investigate the $\lambda$-NRW on the augmented trees; we identify the underlying set with the hyperbolic boundary and the Martin boundary. The estimations of the Martin kernel and the Na\"im kernel are stated. Finally, in Section 5 we discuss the induced Dirichlet forms and a criterion of the critical exponents.

\bigskip

\section{Preliminaries}
\label{sec:2}

Let $(X,E)$ be a countably infinite connected graph with an edge set $E$ and a {\it root}  $\vartheta \in X$. For $x \in X$, we let $|x|:=d(\vartheta,x)$ where $d(\cdot, \cdot)$ is the graph distance, and let  $\mathcal J_n =\{x\in X: |x|=n\}$ be the $n$-th level of the graph. We define a partial order $\prec$ on $X$ with $y \prec x$ if $x$ lies on some geodesic path $\pi(\vartheta,y)$;  for $m \geq 0$,  and define the {\it $m$-th descendant set} of $x$ as ${\mathcal J}_m(x):= \{y \in X: y \prec x \hbox{ and } |y|=|x|+m\}$.

\medskip
Let $E_v:=\{(x,y) \in E: |x|-|y|=\pm 1\}$ and $E_h:=\{(x,y) \in E: |x|=|y|\}$ denote the {\it vertical edge} set and the {\it horizontal edge} set respectively. Clearly $E=E_v \cup E_h$. We call $(X,E)$ a {\it (rooted) tree} if $E=E_v$ and any vertex $x \in X\setminus \{\vartheta\}$ has a unique {\it parent} $x^-$ that satisfies $x \in {\mathcal J}_1(x^-)$.

\bigskip

To discuss the hyperbolicity of an infinite rooted graph $(X,E)$, let us recall the setup by Gromov \cite{Gr} (see also \cite{Wo}). The {\it Gromov product} of two vertices $x,y \in X$ is $(x|y) := \frac 12 \big(|x|+ |y| - d(x,y)\big)$.
If there exists $\delta \geq 0$ such that $(x|y) \geq \min \{(x|z),(z|y)\}-\delta$ for all $x,y,z \in X$, then we say that $(X,E)$ is {\it (Gromov) hyperbolic}.
Clearly every tree is hyperbolic with $\delta=0$. For $a>0$, we define the {\it Gromov metric} $\varrho_a (\cdot, \cdot)$ on $X$ by
\begin{equation} \label{eq2.1}
\varrho_a(x, y) = \begin {cases}  e^{-a(x|y)},\quad  & \hbox {if} \ x \neq y,\\ 0 ,  & \hbox {if}\  x=y.\end{cases}
\end{equation}
Then $\varrho_a(x, y) \leq e^{\delta a}\max \{\varrho_a(x, z), \varrho_a(z,y)\}$ (hence $\varrho_a$ is a quasi-metric as defined in \eqref{eq3.0}); this $\varrho_a$ is equivalent to a metric when $e^{\delta a} < \sqrt 2$.

\medskip

\begin{definition} \label{de2.2}
For $a>0$, let $\widehat{X}_H$ denote the $\varrho_a$-completion of $X$, and call $\partial_H X := \widehat{X}_H \setminus X$ the {\it hyperbolic boundary} of $(X,E)$.
\end{definition}

We note that the topology does not depend on the value of $a$.  We treat $\xi \in \partial_H X$  as the limit of a $\rho_a$-Cauchy sequence $\{x_n\}_{n=0}^\infty$. In particular, we can take $\{x_n\}_{n=0}^\infty$ to be a {\it geodesic ray} $[x_0, x_1, \cdots ]$ (i.e., each finite segment is a geodesic of the ambient graph  $(X, E)$). We denote such a ray by $\pi(x_0, \xi)$. We extend the partial order $\prec$ by writing $\xi \prec x$ if $x$ lies on some $\pi(\vartheta,\xi)$,   and it is easy to define geodesics for  $x \in X, \xi \in \partial_H X$, and (bi-infinite) geodesics for $\xi, \eta \in \partial_H X$.

For the augmented tree discussed in the next section, we can extend the Gromov product to  $(\xi | \eta), \ \xi, \eta \in   \partial_H X$ (see \eqref{eq3.4}), then we can use this to define  $\varrho_a(\cdot,\cdot)$ on ${\widehat X}_H$ in the same way as in \eqref{eq2.1}, and it is still a  quasi-metric.

\medskip

In \cite{Ka2}, Kaimanovich introduced a class of rooted graphs $(X,E)$ that we call {\it pre-augmented trees}:  the subgraph $(X,E_v)$ is a tree, and
$$
(x,y) \in E_h \ \Rightarrow \  x^-=y^-  \hbox { or } (x^-,y^-) \in E_h.
 $$
 The pre-augmented tree allows us to have a simple expression of its geodesics.
 For two distinct $x,y$ in $X$, there exists a geodesic consisting of three segments: $\pi (x, u)$, $\pi (u, v)$ and $\pi(v, y)$, where $u,v$  satisfy $x \prec u$, $y \prec v$, $|u| = |v|$ (hence the middle one is horizontal, and the other two are vertical). In this paper we always use geodesics of this form unless otherwise stated. Note that for $x, y \in X$, the geodesic $\pi(x, y)$ is not necessary unique, and we call the one whose horizontal part satisfies $\pi(u, v) \subset \mathcal{J}_n$ with $n$ smallest a {\it canonical geodesic} between $x$ and $y$. 

 \medskip

 The following is a useful criterion for  a pre-augmented tree to be hyperbolic.

\medskip

\begin{proposition} \label{th2.4} \hspace{-1mm} {\rm \cite[Theorem 2.3]{LW1}}
A pre-augmented tree is hyperbolic if and only if there exists $M<\infty$ such that the lengths of all horizontal geodesic segments are bounded by $M$.
\end{proposition}

For a hyperbolic pre-augmented tree, it is direct to show that for $\xi, \eta \in \partial_H X$, there is a canonical geodesic $\pi (\xi, \eta)$.

\bigskip

Next we consider  random walks and  their Martin boundaries.  An {\it (electric) network} $(X, E,c)$ is a locally finite connected graph $(X,E)$ together with a nonnegative function $c$ on $X \times X$ that satisfies $c(x,y)=c(y,x)>0$ if and only if  $(x,y) \in E$; we call such $c(x,y)$ the {\it conductance} of the edge $(x,y)$, and $m(x):= \sum_{y \in X} c(x,y) > 0$ the {\it total conductance} at $x$. Let $\ell(X)$ denote the collection of real functions on $X$. The {\it graph energy} of $f \in \ell(X)$ is defined by
\begin{equation} \label{eq2.3}
\mathcal E_X[f] = \frac 12 \sum_{x,y \in X} c(x,y)|f(x)-f(y)|^2,
\end{equation}
and the domain of $\mathcal E_X$ is $\mathcal D_X = \{f \in \ell(X): \mathcal E_X[f]<\infty\}$. For $X' \subset X$, a function $f \in \ell(X)$ is said to be {\it harmonic} on $X'$ if $\sum_{y \in X}c(x,y)(f(x)-f(y)) = 0$ for any $x \in X'$. We define the {\it effective resistance} between two disjoint finite nonempty subsets $F,G \subset X$ by
\begin{equation} \label{eq2.4}
R_X(F,G) = (\min\{\mathcal E_X[f]: f \in \ell(X) \hbox{ with } f = 1 \hbox{ on } F, \hbox{ and } f=0 \hbox{ on } G\})^{-1}.
\end{equation}
Also set $R_X(F,G)=0$ if $F \cap G \neq \emptyset$ by convention. Clearly the energy minimizer in \eqref{eq2.4} is unique, and is harmonic on $X \setminus (F \cup G)$.

\medskip

Every conductance $c$ defines a {\it reversible random walk} $\{Z_n\}_{n=0}^\infty$ on the graph $(X,E)$, in which the transition probability is given by $P(x,y) = \frac{c(x,y)}{m(x)}$, $x,y \in X$.
In this study we always fix $\vartheta \in X$ so that $(X,E)$ is a rooted graph, and assume that $\{Z_n\}$ is {\it transient}, i.e., its {\it Green function} $G(x,y):=\sum_{n=0}^\infty P^n(x,y) < \infty$  for all $x,y \in X$. Write $\mathbb P(\cdot \mid Z_0=x)=\mathbb P_x(\cdot)$ for short. The {\it ever-visiting probability} is defined as $F(x,y):=\mathbb P_x\{\exists\ n \geq 0 \hbox{ such that } Z_n=y\}$; it satisfies $G(x,y) = F(x,y)G(y,y)$. We define the {\it Martin kernel} by
\begin{equation*}
K(x,y) = \frac{G(x,y)}{G(\vartheta,y)} = \frac{F(x,y)}{F(\vartheta,y)}, \qquad x,y \in X.
\end{equation*}

\begin{definition} \label{de2.5}
We call $\mathcal M = \widehat{X} \setminus X$ the {\it Martin boundary} of $\{Z_n\}$, where $\widehat{X}$ is the  compactification of $X$ such that $K(x,\cdot)$ extends continuously to $\widehat{X}$ for all $x \in X$.
\end{definition}

Under the Martin topology, the random walk $\{Z_n\}$ converges to an $\mathcal M$-valued random variable $Z_\infty$ almost surely. The {\it hitting distribution} $\nu$ is defined by $\nu(B)=\mathbb P_\vartheta(Z_\infty \in B)$ for any Borel set $B \subset \mathcal M$. For $\xi \in \mathcal M$, define the {\it $\xi$-process} to be the random walk $\{Z_n^\xi\}_{n=0}^\infty$ on $(X,E)$ with the transition probability $P^\xi(x,y):=P(x,y)K(y,\xi)/K(x,\xi)$; here $P^\xi$ is stochastic since $K(\cdot,\xi)$ is $P$-harmonic. Analogously we define the hitting distribution  $\nu^\xi$ for the $\xi$-process. The {\it minimal Martin boundary} of $\{Z_n\}$ is given by $\mathcal M_{\min} := \{\xi \in \mathcal M: \nu^\xi \hbox{ is the point mass at } \xi\}$, and satisfies $\nu (\mathcal M\setminus \mathcal M_{\min})=0$.

\bigskip

For $F \subset X$, let $m(F):= \sum_{x \in F}m(x)$ and let $c(\partial F):=\sum_{x \in F,y \notin F} c(x,y)$.
We say that a conductance $c$ on $(X,E)$ has {\it strong isoperimetry} \eqref{eq2.2} if
\begin{equation} \label{eq2.2}
\sup\left\{\frac{m(F)}{c(\partial F)}: F \subset X \hbox{ and } F \hbox{ is finite}\right\} < \infty. \tag{SI}
\end{equation}
Note that \eqref{eq2.2} implies the transience of $\{Z_n\}$, but the converse is not true. As an alternative version of Ancona's Theorem in \cite{An}, we have

\medskip

\begin{theorem} \label{th2.6} \hspace{-1mm} {\rm (Ancona)}
Let $(X,E)$ be a hyperbolic graph. Suppose a conductance $c$ on $(X,E)$ satisfies hypotheses \eqref{eq2.2} and $\inf_{(x,y) \in E}P(x,y) >0$. Then the reversible random walk $\{Z_n\}_{n=0}^\infty$ satisfies
\begin{equation} \label{eq2.5}
F(x,y) \asymp F(x,z)F(z,y)
\end{equation}
uniformly for all $x, y, z \in X$ such that $z$ lies on some $\pi(x, y)$. Moreover, the Martin boundary ${\mathcal M}$  equals the minimal Martin boundary $\mathcal M_{\min}$, and is homeomorphic to the hyperbolic boundary $\partial_HX$.
\end{theorem}



We define the {\it Poisson integral} $H: L^1(\mathcal M, \nu) \to \ell(X)$ by
\begin{equation} \label{eq2.6}
(Hu)(x) = \int_{\mathcal M} K(x,\xi)u(\xi)d\nu(\xi), \qquad u \in L^1(\mathcal M, \nu),\ x \in X.
\end{equation}
Since $K(\cdot,\xi)$ is harmonic for all $\xi \in \mathcal M$, we know that $Hu$ is harmonic on $X$. Via the operator $H$, the energy form $(\mathcal E_X, \mathcal D_X)$ in \eqref{eq2.3} induces a bilinear form $(\mathcal E_{\mathcal M}, \mathcal D_{\mathcal M})$ on $L^2(\mathcal M, \nu)$ defined by
\begin{equation} \label{eq2.7}
\mathcal E_{\mathcal M}(u,v) = \mathcal E_X(Hu,Hv), \qquad u,v \in \mathcal D_{\mathcal M},
\end{equation}
where the domain $\mathcal D_{\mathcal M} = \{u \in L^2(\mathcal M, \nu): Hu \in \mathcal D_X\}$. Define the {\it Na{\"i}m kernel} by
\begin{equation} \label {eq2.7'}
\Theta(x,y) = \frac{K(x,y)}{G(x,\vartheta)} = \frac{F(x,y)}{F(x,\vartheta)G(\vartheta,\vartheta)F(\vartheta,y)}, \qquad x,y \in X.
\end{equation}
This kernel is symmetric on $X \times X$, and can be extended continuously to $\Theta(x,\eta)$ on $X \times \mathcal M$ as the Martin kernel $K(x,y)$ does. The extension to $\Theta( \xi, \eta)$ for $\xi, \eta \in {\mathcal M}$ is more involved. With the notations in the rooted graph $(X,E)$, for $m \geq 0$ and $z \in X$, let $\tau_m^*$ be the last visit time of $\mathcal J_m$ by the $\xi$-process $\{Z_n^\xi\}$, i.e., $\tau_m^*=\sup\{n \geq 0: |Z_n^\xi|=m\}$, and let $\ell_m^\xi(z) = \mathbb P_\vartheta^\xi(Z_{\tau_m^*}^\xi = z)$. The following extension of $\Theta(\cdot,\cdot)$ on $\mathcal M \times \mathcal M \setminus \Delta$ (here $\Delta:=\{(\xi,\xi): \xi \in \mathcal M\}$) was introduced by Silverstein in \cite{Si}:
\begin{equation} \label{eq2.8}
\Theta(\xi,\eta) = \lim_{m \to \infty} {\sum}_{z \in \mathcal J_m} \ell_m^\xi(z)\Theta(z,\eta), \qquad \xi \neq \eta \in \mathcal M.
\end{equation}
This limit exists since the sum is increasing in $m$. One of the main results in \cite{Si} is stated below.

\begin{theorem} \label{th2.7} \hspace{-1mm} {\rm (Silverstein)}
The induced form $(\mathcal E_{\mathcal M}, \mathcal D_{\mathcal M})$ in \eqref{eq2.7} has the expression
\begin{equation*}
\mathcal E_{\mathcal M}(u,v) = \frac{m(\vartheta)}{2} \iint_{\mathcal M \times \mathcal M \setminus \Delta} (u(\xi)-u(\eta))(v(\xi)-v(\eta))\Theta(\xi,\eta)d\nu(\xi)d\nu(\eta), \quad u,v \in \mathcal D_{\mathcal M}.
\end{equation*}
Moreover, $\mathcal D_{\mathcal M} = \{u \in L^2(\mathcal M, \nu): \mathcal E_{\mathcal M}[u]:= \mathcal E_{\mathcal M}(u,u)<\infty\}$.
\end{theorem}

\bigskip

\section{Partition and augmented tree}
\label{sec:3}


The aim of this section is to associate every compact space of homogeneous type $(K,\rho,\mu)$ with some hyperbolic graph $(X,E)$ such that the hyperbolic boundary $(\partial_HX,\varrho_a) \approx (K,\rho)$.
%
 Recall that a {\it quasi-metric} $\rho$ on a set $M$ is a function on $M\times M$ to $[0, \infty)$ that satisfies the conditions of a metric except that the triangle inequality is weakened  \cite {AM}: for some $C_\rho \geq 1$,
\begin{equation} \label{eq3.0}
\rho(\xi , \eta) \leq C_\rho \big (\rho (\xi, \zeta) + \rho (\zeta, \eta)\big), \quad \forall \ \xi , \eta, \zeta \in M.
\end{equation}
 A quasi-metric $\rho$ defines a topology $\mathcal T_\rho$ on $M$ by
$$
\mathcal T_\rho = \{\Omega \subset M: \forall \, \xi \in \Omega,\, \exists\,\varepsilon>0 \ {\rm s.t.} \ B(\xi,\varepsilon):=\{\eta \in M: \rho(\xi,\eta)<\varepsilon\} \subset \Omega\}.
$$
Note that for some $\rho$, not all balls $B(\xi,r)$ are open with respect to $\mathcal T_\rho$. However, Macias and Segovia \cite{MS} proved that there exists a quasi-metric $\rho'$, equivalent to $\rho$ in the sense that $\rho'(\xi,\eta) \asymp \rho(\xi,\eta)$ (hence $\mathcal T_{\rho'} = \mathcal T_\rho$), such that all associated metric balls belong to $\mathcal T_\rho$. Without loss of generality we will assume that metric balls are open throughout the paper.

\medskip

\begin{definition} \label{de3.1}
A space of homogeneous type is a set $M$ equipped with a quasi-metric $\rho$  for which all balls $B(\xi,r)$ are open, and a non-negative regular Borel measure $\mu$ that satisfies the
doubling condition: there exists $C>0$ such that
\begin{equation} \label{eq3.1}
0 < \mu(B(x,2r)) \leq C  \mu(B(x,r)) < \infty, \quad \forall\ \xi \in M, \, r>0. \tag{VD}
\end{equation}
\end{definition}

Typical examples of spaces of homogeneous type include classical Euclidean domains, Riemannian manifolds and homogeneous spaces. A non-isotropic example is $M= {\mathbb R}^d$ with the Lebesgue measure $\mu$ and $\rho(x, y) =\sum_{i=1}^d |x_i-y_i|^{s_i}$, where $s_1,\cdots, s_d$ are positive numbers.
For more examples and a general method for constructing these spaces, the reader can refer to \cite{CW}.
In fractal geometry, an important class of examples is given by the $\alpha$-sets, which are closed subsets in ${\mathbb R}^d$ equipped with the Euclidean distance and the Ahlfors $\alpha$-regular measures $\mu$, i.e., $\mu (B(x,r)) \asymp r^\alpha$ for all $x \in K$  and $r \in (0,1)$.

\bigskip

Through out this paper we will use $(K, \rho, \mu)$ to denote a compact space of homogeneous type unless otherwise specified. For convenience, we assume that $\diam(K) := \sup_{x, y \in K} \rho(x, y) =1$ and $\mu$ is a probability measure. We first introduce a {\it tree of partitions} on the space as in \cite{LW2}.

\begin{definition} \label{de3.2}
Let $\mathcal K$ denote the collection of nonempty compact subsets of $K$. We call a triple $(X,E_v,\Phi)$ an {\it index tree} of $(K,\rho,\mu)$ if $(X,E_v)$ is a  tree together with a set-valued mapping $\Phi: X \to \mathcal K$ that satisfies the following:
\begin{enumerate}[(i)]
\item $\Phi(\vartheta)=K$, and $\Phi(x) = \bigcup_{y \in \mathcal J_1(x)} \Phi(y)$ for all $x \in X$;

\item $\exists \, b, c>0$ and $r_0 \in (0,1)$ such that each $x \in X$ satisfies $B(\xi_x, b r_0^{|x|}) \subset \Phi(x) \subset B(\xi_x,c r_0^{|x|})$ for some $\xi_x \in K$;

\item ({\it $\mu$-separation}) for distinct $x,y \in X$ with $|x|=|y|$, $\mu(\Phi(x) \cap \Phi(y)) = 0$.
\end{enumerate}
\end{definition}

The parameter $r_0$ in (ii) is referred to as the {\it contraction ratio} of $\Phi$. Note that for each $\xi \in K$, there exists a sequence of $x_n \in {\mathcal J}_n$ such that $\xi\in \Phi(x_n)$ for all $n$, which yields $\{\xi\} \in \bigcap_{n=0}^\infty \Phi (x_n)$, and we call $\{x_n\}_{n=0}^\infty$ a {\it representation} of $\xi \in K$.

\medskip

We first show that a compact space of homogeneous type always admits an index tree according to the above definition. This follows from the construction of {\it dyadic cubes} by Christ \cite{Ch}.

\medskip

\begin {example} \label {ex3.1} {\rm {\bf Dyadic cubes}\  For $\varepsilon>0$, we define an {\it $\varepsilon$-net} on $K$ as a finite subset $\Xi \subset K$ such that $\bigcup_{\xi \in \Xi} B(\xi,\varepsilon) = K$, and $\rho(\xi,\eta) \geq \varepsilon$ if $\xi,\eta \in \Xi$ with $\xi \neq \eta$.
Assume that $K$ has diameter 1, and fix $r_0 \in (0,1)$ and $b > 0$ such that \ $\frac {C_\rho^3 r_0} {1-C_\rho r_0} + C_\rho^2 b \leq \frac 12$.
For $m \geq 0$, let $\Xi_m$ be an $r_0^m$-net on $K$, and index its elements by $\mathcal J_m$ so chosen that the $\mathcal J_m$'s are disjoint;  then we write $\Xi_m=\{\xi_x\}_{x \in \mathcal J_m}$.

\vspace{1mm}

For $m \geq 1$ and $x \in \mathcal J_m$, we choose $x^- \in \mathcal J_{m-1}$ such that $\xi_{x^-}$ is the nearest point from $\xi_x$ among all points in $\Xi_{m-1}$ (if there are two or more nearest points, we select an arbitrary one from them); clearly $\rho(\xi_x,\xi_{x^-}) \leq r_0^{m-1}$. Let $X = \bigcup_{m=0}^\infty \mathcal J_m$ and let $E_v=\{(x,x^-),(x^-,x):x \in X, \, x \neq \vartheta\}$. Then $(X,E_v)$ is a  tree that satisfies $\mathcal J_1(x) \neq \emptyset$ for every $x \in X$. 
We define
\begin{equation*}
\Phi^\circ(x) = {\bigcup}_{z: z \prec x} B(\xi_z, b r_0^{|z|}) \quad \hbox{and} \quad \Phi(x) = \overline{\Phi^\circ(x)}, \qquad x \in X,
\end{equation*}
and call $\{\Phi(x)\}_{x \in X}$ a set of {\it dyadic cubes}, where $\prec$ is the partial order defined in the last section. Then it is clear that $K=\overline{\Phi^\circ(\vartheta)}=\Phi(\vartheta)$,   $\Phi^\circ(x) \subset B(\xi_x,\frac {C_\rho r_0^{|x|}}{1- C_\rho r_0})$ for all $x \in X$, and $\Phi(x) = \bigcup_{y \in \mathcal J_1(x)} \Phi(y)$ for all $x \in X$.

\vspace {0.1cm}
Using \eqref{eq3.1} and the Lebesgue differentiation theorem, it was proved in \cite {Ch} that $\mu(\bigcup_{x \in \mathcal J_m} \Phi^\circ(x)) = \mu(K)$ for all $m$; moreover, for $x \neq y \in X$ with $|x|=|y|$, it follows from
 $\frac {C_\rho^3 r_0} {1-C_\rho r_0} + C_\rho^2 b \leq \frac 12$ that $\Phi^\circ(x) \cap \Phi^\circ(y) = \emptyset$, which implies $\mu(\Phi(x) \cap \Phi(y))=0$.
}\end{example}

\bigskip
 The original framework of index trees in Definition \ref{de3.2}, as we will see below, is the natural tree structures on self-similar sets, which are known as the symbolic representations of their associated iterated function systems.

\medskip

\begin{example} \label{ex3.3} \hspace{-2mm} {\bf Self-similar sets} {\rm
Let $\{S_i\}_{i \in \Sigma}$, where $\Sigma$ is a finite set, be an iterated function system (IFS) of contractive similitudes on $\mathbb R^d$ such that $S_i$ has contraction ratio $r_i \in (0,1)$. The {\it self-similar set} is the unique nonempty compact set $K \subset \mathbb R^d$ that satisfies $K = \bigcup_{i \in \Sigma} S_i(K)$. We assume that the IFS satisfies the open set condition (OSC), i.e., there exists a nonempty bounded open set $O$ such that $S_i(O) \subset O$, and $S_i(O) \cap S_j(O) = \emptyset$ for $i \neq j \in \Sigma$. Let $\rho$ be the Euclidean metric and $\mu$ be the $\alpha$-Hausdorff measure on $K$, where $\alpha$ is the Hausdorff dimension of $K$ (determined by $\sum_{i \in \Sigma} r_i^\alpha = 1$).

\vspace{1mm}

Let $\Sigma^*=\bigcup_{n \geq 0} \Sigma^n$ (by convention $\Sigma^0 =\{\vartheta\}$). For $x=i_1 i_2 \cdots i_k \in \Sigma^*$, write $S_x := S_{i_1} \circ S_{i_2} \circ \cdots \circ S_{i_k}$ and $r_x := r_{i_1} r_{i_2} \cdots r_{i_k}$ for short. Let $r_0:=\min_{i \in \Sigma} r_i$. Define $\mathcal J_0 = \{\vartheta\}$ and for $n \geq 1$,
\begin{equation*}
\mathcal J_m = \{x=i_1i_2 \cdots i_k \in \Sigma^*: r_x \leq r_0^m < r_{i_1i_2 \cdots i_{k-1}}\}.
\end{equation*}
For $m \geq 1$ and $x \in \mathcal J_m$, there exists a unique $x^- \in \mathcal J_{m-1}$ such that $S_x(K) \subset S_{x^-}(K)$. Naturally we set $X = \bigcup_{m=0}^\infty \mathcal J_m$, $E_v = \{(x,x^-),(x^-,x): x \in X, \, x \neq \vartheta\}$ and $\Phi(x) = S_x(K)$ for $x \in X$. It is straightforward that $(X,E_v,\Phi)$ is an index tree of $(K,\rho,\mu)$ with contraction ratio $r_0$.
}
\end{example}

\noindent {\bf Remark.} For some self-similar sets we can replace the Hausdorff measure by a {\it self-similar measure} $\mu$ defined by $\mu(\cdot) = \sum_{i \in \Sigma} p_i\mu(S_i^{-1}(\cdot))$, where $\{p_i\}_{i \in \Sigma}$ is a set of positive probability weights such that $\mu$ is a doubling measure on $K$ (see \cite{Y}).

\medskip
%
%
%


An index tree provides a convenient symbolic representation of elements in the underlying space $K$. However, this representation does not reflect the actual geometric structure of $K$. For this reason, we strengthen the index tree by adding a set of horizontal edges in  light of the positional relationship among the subsets in each $\Phi(\mathcal J_n)$ (see \cite{Ka2,LW1,LW2}).

\medskip

\begin{definition} \label{de3.4}
Let $(X,E_v,\Phi)$ be an index tree of $(K,\rho,\mu)$. For a fixed $\gamma>0$, define a horizontal edge set on $X$ by
\begin{equation} \label{eq3.3}
E_h := {\bigcup}_{m=1}^\infty \{(x,y) \in \mathcal J_m \times \mathcal J_m:  x \neq y,\ {\rm dist} (\Phi(x),\Phi(y)) \leq \gamma \cdot r_0^m\}.
\end{equation}
Let $E = E_v \cup E_h$. We call $(X,E,\Phi)$ an {\it augmented tree} of $(K,\rho,\mu)$.
\end{definition}

\noindent {\bf Remark.} The notion of augmented tree was coined  by Kaimanovich in \cite{Ka2} where the inequality in \eqref{eq3.3} was originally given by $\Phi(x) \cap \Phi(y) \neq \emptyset$. The modified version we define here is adopted from \cite{LW2} in order to avoid the superfluous conditions in \cite{LW1}.

\bigskip

Clearly $(X,E)$ is a pre-augmented tree, and  the hyperbolicity can be verified using Proposition \ref{th2.4} and a similar proof as in \cite[Theorem 1.2]{LW2}. On the augmented tree, We extend the Gromov product to the boundary $\partial_HX$ by
\begin{equation} \label{eq3.4}
(x|\eta) = \sup\left\{\lim_{n \to \infty}(x|y_n)\right\}, \ \ (\xi|\eta) = \sup\left\{\lim_{n \to \infty}(x_n|y_n)\right\}, \quad x \in X,\ \xi,\eta \in \partial_HX,
\end{equation}
where the supremum is taken over all geodesic rays $\pi(\vartheta,\xi) = [\vartheta, x_1,\cdots]$ and $\pi(\vartheta,\eta) = [\vartheta, y_1,\cdots]$; the above limits exist since both $(x|y_n)$ and $(x_n|y_n)$ are increasing in $n$, and the difference of the limits for different rays are at most $1$ (see \cite{LW2,KLW} for details).
 Then $\varrho_a(\cdot,\cdot)$ on ${\widehat X}_H$ can be defined in the same way as in \eqref{eq2.1}, and it is still a quasi-metric.

\medskip

We define a map $\kappa$ on the collection of geodesic rays starting at $\vartheta$ to $K$ by
\begin{equation} \label{eq3.2}
\{\kappa([x_0,x_1,\cdots])\} = {\bigcap}_{n=0}^\infty \Phi(x_n).
\end{equation}
For each $\xi \in \partial_HX$, it is known that $\kappa(\pi(\vartheta,\xi))$ is independent of the choice of $\pi(\vartheta,\xi)$, and hence $\kappa$ defines a map $\hat{\kappa}$ from $\partial_HX$ to $K$ by $\hat{\kappa}(\cdot) = \kappa(\pi(\vartheta,\cdot))$, which has been proved to be a H{\"o}lder continuous homeomorphism (see \cite[Theorem 1.3]{LW2}). The following is the main reason for introducing the augmented tree.

\medskip

\begin{theorem} \label{th3.5} \hspace{-2mm}
Let $(X,E,\Phi)$ be an augmented tree of $(K,\rho,\mu)$. Then $(X,E)$ is hyperbolic, and $(K,\rho)$ is H{\"o}lder equivalent to the hyperbolic boundary $(\partial_HX,\varrho_a)$ in the sense that $\varrho_a(\xi,\eta) \asymp \rho(\hat\kappa(\xi),\hat\kappa(\eta))^{-a/\log r_0}$ for all $\xi,\eta \in \partial_HX$, where $r_0$ is the contraction ratio of $\Phi$.
\end{theorem}
\noindent {\bf Remark.} This theorem is still valid even if we omit the condition ``$\mu$-separation" in Definition \ref{de3.2} (see \cite {LW2}).

\bigskip

From this result we can identify $\partial_HX$ with $K$. By using the doubling condition of $\mu$ and the triangle inequality of the quasi-metric, we have two elementary results on $\mu(\Phi(x))$ which will play an important role in the study of random walks on the augmented trees.

\medskip

\begin{proposition} \label{th3.6}
Let $(X,E,\Phi)$ be an augmented tree of $(K,\rho,\mu)$. Then there is a constant $C>1$ such that
\begin{equation} \label{eq3.5}
C^{-d(x,y)} \mu(\Phi(x)) \leq \mu(\Phi(y)) \leq C^{d(x,y)} \mu(\Phi(x)), \qquad \forall \ x,y \in X.
\end{equation}
Consequently every augmented tree has bounded degree.
\end{proposition}

\medskip

%
%

\begin{proposition} \label{th3.8}
Let $(X,E,\Phi)$ be an augmented tree of $(K,\rho,\mu)$. Define for  $\xi \neq \eta \in X \cup K$,
\begin{equation} \label{eq3.5'}
p_\mu(\xi,\eta) := \sup \{\mu(\Phi(z)): z \in X \hbox{ and } z \hbox{ lies on some canonical } \pi(\xi,\eta)\}.
\end{equation}
Then the supremum can be reached at a vertex lying on the horizontal segment of some canonical geodesic between $\xi$ and $\eta$, and there exists $C > 1$ such that
\begin{equation*} \label{eq3.5''}
p_\mu(\xi,\eta) \leq C \max\{p_\mu(\xi,\sigma),p_\mu(\sigma,\eta)\}
\end{equation*}
for all different $\xi,\eta,\sigma \in X \cup K$.
Moreover,   $p_\mu(\cdot,\cdot)$ satisfies the estimate
\begin{equation*} \label{eq3.6}
p_\mu(\xi,\eta) \asymp V(\xi,\eta) := \mu(B(\xi,\rho(\xi,\eta))), \qquad \forall\ \xi,\eta \in K,\ \xi \neq \eta.
\end{equation*}
\end{proposition}

\medskip

\section{Reversible random walks on augmented trees}
\label{sec:4}

 Let $(X,E,c)$ be a (rooted) network with energy $\mathcal E_X[f] = \frac 12 \sum_{x,y \in E} c(x,y)|f(x)-f(y)|^2, \ f\in \ell(X)$; it induces a reversible random walk $\{Z_n\}_{n=0}^\infty$. For each $x \in X \setminus \{\vartheta\}$,  we define the {\it return ratio} at $x$ by
\begin{equation} \label{eq4.1}
\lambda(x) := \frac{\mathbb P_x(|Z_1| = |x|-1)}{\mathbb P_x(|Z_1| = |x|+1)} = \frac{\sum_{z: x \in \mathcal J_1(z)}c(x,z)}{\sum_{y: y \in \mathcal J_1(x)}c(x,y)}.
\end{equation}

\medskip

\begin{theorem} \label{th4.1}
Let $(X,E, c)$ be a network, and suppose the associated $\{Z_n\}_{n=0}^\infty$ satisfies \\
(i) $  {\inf}_{(x,y) \in E}\, P(x,y) > 0$; \  (ii)  $ {\sup}_{x \in X \setminus \{\vartheta\}} \, \lambda(x) < 1$. Then the conductance on $(X,E$) has strong isoperimetry \eqref{eq2.2}.

\vspace {0.1cm}
In particular, if the above $\{Z_n\}$ is defined on an augmented tree $(X,E,\Phi)$ and satisfies (i) and (ii), then the hypotheses in Theorem \ref{th2.6} are fulfilled. Hence the Martin boundary $\mathcal M=\mathcal M_{\min}$, and is homeomorphic to the hyperbolic boundary $\partial_HX$ as well as the underlying compact set $K$.
\end{theorem}

The proof of the (SI) is similar to \cite[Theorem 5.1]{KLW}. Henceforth, we will identify $K$ with the Martin boundary $\mathcal M$.
%
To consider the hitting distribution $\nu$, we introduce a key condition on $\{Z_n\}_{n=0}^\infty$ which is satisfied by the $\lambda$-NRW introduced below:
\begin{equation} \label{eq4.4}
\mathbb P_\vartheta(Z_{\tau_m} = x) = \mu(\Phi(x)), \qquad \forall \ m \geq 1, \ x \in \mathcal J_m, \tag{HD}
\end{equation}
where $\tau_m:= \inf\{n \geq 0: |Z_n|=m\}$ is the first hitting time of $\mathcal J_m$ by $\{Z_n\}_{n=0}^\infty$.

\medskip

\begin{proposition} \label{th4.2}
In addition to the  assumptions in Theorem \ref{th4.1}, suppose (HD) is satisfied. Then the hitting distribution $\nu$ equals $\mu$.
\end{proposition}

One of our objectives is to obtain some explicit estimates of the  Martin kernel and the Na{\"i}m kernel. For this we consider a class of reversible random walks with {\it constant return ratio} $\lambda$:
\begin{equation} \label{eq4.5}
\lambda(x) \equiv \lambda, \qquad \forall\ x \in X \setminus \{\vartheta\}, \tag{R$_\lambda$}
\end{equation}
where $\lambda \in (0,1)$ is a constant. With this hypothesis, a subsequence of $\{|Z_n|\}_{n=0}^\infty$ is seen to be a birth and death chain by counting the time instants only when $\{Z_n\}$ moves upward or downward (see \cite[Proposition 4.1]{KLW}). With this in mind we introduce the following key definition.

\medskip

\begin{definition} \label{de4.3}
Let $(X,E,\Phi)$ be an augmented tree of $(K,\rho,\mu)$. A reversible random walk $\{Z_n\}$ on $(X,E)$ is said to be {\rm $\lambda$-natural} ($\lambda$-NRW) if it satisfies
\begin{equation*} \label{eq4.6}
c(x,x^-) = \lambda^{-|x|}\mu(\Phi(x)), \quad \forall \ x \in X \setminus \{\vartheta\}, \quad \hbox{and} \quad c(x,y) \asymp c(x,x^-), \quad \forall \ (x,y) \in E_h.
\end{equation*}
\end{definition}

The notion of natural random walks was first introduced and studied in \cite{KLW} for the self-similar sets in Example \ref{ex3.3}. In that case the $\alpha$-Hausdorff measure $\mu$ is the self-similar measure with the ``natural weight" $p_i = r_i^\alpha$ where the $r_i$'s are the contraction ratios of the corresponding IFS, and  $\mu(\Phi(x)) = p_{i_1}\cdots p_{i_n}$ for $x= i_1\cdots i_n$. Also in \cite {KLW}, we defined a class of quasi-NRW with $\mu$ a self-similar measure, and necessarily
$\mu$ is doubling \cite [Theorem 4.8] {KLW}. Hence the quasi-NRW is also a NRW in Definition \ref{de4.3}.

\bigskip
It is easy to show that for any  $\lambda$-NRW the condition $(R_\lambda)$ is satisfied. Also the bounded degree property of the augmented tree (Proposition \ref{th3.6}) implies that every augmented tree admits a $\lambda$-NRW, and $ \inf_{(x, y)\in E} P(x, y) = \inf_{(x, y)\in E} \frac{c(x, y)}{m(x)} >0$.  In view of Theorem \ref{th4.1} and Proposition \ref{th4.2}, we have

\medskip

\begin{theorem} \label{th4.4} \hspace{-2mm}
Let $(X,E,\Phi)$ be an augmented tree of $(K,\rho,\mu)$ and let $\{Z_n\}_{n=0}^\infty$ be a $\lambda$-NRW on $(X,E)$. Then the Martin boundary, the hyperbolic boundary and $K$ are homeomorphic.

\vspace {0.1cm}
Moreover $\{Z_n\}_{n=0}^\infty$ satisfies property (HD), and
\begin{equation*} \label{eq4.7}
F(x,\vartheta) = \lambda^{|x|}, \quad \hbox{and} \quad F(\vartheta,x) \asymp \mu(\Phi(x)), \qquad \forall \ x \in X,
\end{equation*}
and the hitting distribution $\nu$ on $K$ is $\mu$.
\end{theorem}

\medskip

\begin{theorem} \label{th4.5}
Let $\{Z_n\}_{n=0}^\infty$ be a $\lambda$-NRW on an augmented tree $(X,E,\Phi)$. Then the ever-visiting probability obeys the estimate
\begin{equation*} \label{eq4.8}
F(x,y) \asymp \lambda^{|x|-(x|y)} \mu(\Phi(y))p_\mu(x,y)^{-1}, \qquad\forall \ x,y \in X,
\end{equation*}
where $p_\mu(\cdot,\cdot)$ is defined by \eqref{eq3.5'}. Consequently the Martin kernel satisfies
\begin{equation} \label{eq4.9}
K(x,\eta) \asymp \lambda^{|x|-(x|\eta)}p_\mu(x,\eta)^{-1},  \quad\forall \ x \in X, \ \eta \in X \cup K.
\end{equation}
\end{theorem}


For a natural random walk, as $K \approx \mathcal M$ and $\mu=\nu$ (Theorem \ref{th4.4}), we view the Poisson integral $H$ in \eqref{eq2.6} as an operator on all $\mu$-integrable functions on $K$. Let $C(K)$ denote the collection of continuous functions on $K$. Applying the estimate  in \eqref{eq4.9}, we obtain the same Fatou-type result as in \cite[Corollary 3.2]{KL}.

\medskip
\begin{corollary} \label{th4.5'}
Let $\{Z_n\}_{n=0}^\infty$ be a $\lambda$-NRW on an augmented tree $(X,E,\Phi)$. Then for $u \in C(K)$ and $\varepsilon>0$, there exists a positive integer $n_0$ such that for any $\xi \in K$,
\begin{equation*}
|Hu(x)-u(\xi)| \leq \varepsilon, \qquad \forall\ x \in {\bigcup}_{m \geq n_0} {\mathcal J_m}  \ \hbox { with } \ \xi \in \Phi(x).
\end{equation*}
\end{corollary}

\medskip

It is not difficult to obtain  $\Theta(x,\eta) \asymp \lambda^{-(x|\eta)}p_\mu(x,\eta)^{-1}$ for $x \in X, \ \eta  \in K$ from the above and \eqref{eq2.7'}. However the extension to $x\in X$ to $\xi \in K$  requires more delicate work. By analyzing the limit in \eqref{eq2.8} as in \cite[Theorem 6.3]{KLW}, we can extend the Na{\"i}m kernel estimates to $(K \times K) \setminus \Delta$.

\medskip

\begin{theorem} \label{th4.6}
Let $(X,E,\Phi)$ be an augmented tree of $(K,\rho,\mu)$ with contraction ratio $r_0$, and let $\{Z_n\}_{n=0}^\infty$ be a $\lambda$-NRW on $(X,E)$. Then the Na{\"i}m kernel obeys the estimate
\begin{equation*}
\Theta(\xi,\eta) \asymp \lambda^{-(\xi|\eta)}p_\mu(\xi,\eta)^{-1}, \qquad \forall \ \xi, \eta \in K, \ \xi \neq \eta.
\end{equation*}
Consequently, by Theorem \ref{th3.5} and Proposition \ref{th3.8}, we have
\begin{equation*}
\Theta(\xi,\eta) \asymp \frac{1}{V(\xi,\eta)\rho(\xi,\eta)^\beta}, \qquad \forall \ \xi, \eta \in K, \ \xi \neq \eta,
\end{equation*}
where $V(\xi,\eta) := \mu(B(\xi,\rho(\xi,\eta)))$, and $\beta = \log \lambda / \log r_0$.
\end{theorem}

\medskip

\section{Induced Dirichlet forms}
\label{sec:5}

%

We first recall the definition of a Dirichlet form $\mathcal E$ on $L^2(M, \mu)$, and the regularity of $\mathcal E$ where $M$ is equipped with a quasi-metric. The reader can refer to  \cite{CF, FOT} for the theory of Dirichlet forms.

\begin{definition} \label{de5.1}  Let  $\mu$ be a positive $\sigma$-finite regular Borel measure on $M$ such that  {\rm supp}$(\mu)= M$.
  A {\it Dirichlet form}  $\mathcal E$ with domain ${\mathcal F}$ is a symmetric bilinear form  which is non-negative definite, closed, densely defined on $L^2(M, \mu)$, and satisfies the Markovian property: $u\in {\mathcal F} \Rightarrow \tilde u := (u\vee 0)\wedge 1 \in \mathcal F$ and ${\mathcal E} [\tilde u] \leq {\mathcal E} [u]$. (Here  $\mathcal E [u] := {\mathcal E}(u, u)$ denote the {\it energy} of $u$.)

\vspace{0.1cm}
If further $M$ is given a  locally compact separable quasi-metric $\rho$, then we say that a Dirichlet form $\mathcal E$ is {\it regular} if ${\mathcal F}\cap C_0(M)$ is dense in $C_0(M)$ with the supremum norm, and dense in ${\mathcal F}$ with the ${\mathcal E}^{1/2}_1$-norm.  It is called {\it local} if ${\mathcal E}(u, \upsilon) =0$ for $u, \upsilon \in {\mathcal F}$  having  disjoint compact supports.
\end{definition}

\medskip

For brevity, we write $L^2 = L^2(K,\mu)$, and use $\fint_B$ to denote $\frac{1}{\mu(B)}\int_B$ for any measurable set $B \subset K$ with $\mu(B) > 0$. For $s>0$, we define
\begin{equation} \label{eq5.2}
{\mathcal N}^s_{2,2}[u] := \int_0^\infty \frac{dr}{r} r^{-2s} \int_K \fint_{B(\xi,r)} |u(\xi)-u(\eta)|^2 d\mu(\eta) d\mu(\xi), \qquad u \in L^2(K, \mu),
\end{equation}
and the {\it Besov space} $\Lambda_{2,2}^s = \Lambda_{2,2}^s (K,\rho,\mu) = \{u \in L^2: {\mathcal N}_{2,2}^s(u) < \infty\}$ with the associated norm $\Vert u\Vert_{\Lambda_{2,2}^s} = \Vert u\Vert_{L^2} + ({\mathcal N}_{2,2}^s[u])^{1/2}$. Clearly $\Lambda_{2,2}^{s+\varepsilon} \subset \Lambda_{2,2}^s$ for $\varepsilon>0$, and $\Lambda_{2,2}^s$ can be trivial when $s$ is a large value (e.g. for $K = [0,1]^d \subset \mathbb R^d$, $\Lambda_{2,2}^1$ contains only constant functions).
As a consequence of Silverstein's Theorem (Theorem \ref{th2.7}), Theorems \ref{th4.1}, \ref{th4.2} and \ref{th4.6}, we have

\medskip

\begin{theorem} \label{th5.2}
Let $(K,\rho,\mu)$ be a compact space of homogeneous type, and let $\{Z_n\}$ be a $\lambda$-NRW on the augmented tree $(X,E,\Phi)$. Then the induced form $(\mathcal E_K, \mathcal D_K)$ satisfies
\begin{equation*}
\mathcal E_K[u]:=\mathcal E_X[Hu] \asymp \iint_{K \times K \setminus \Delta} \frac{|u(\xi)-u(\eta)|^2}{V(\xi,\eta)\rho(\xi,\eta)^\beta} d\mu(\xi) d\mu(\eta) \asymp \mathcal N_{2,2}^{\beta/2}[u], \qquad \forall \ u \in \mathcal D_K,
\end{equation*}
and $\mathcal D_K:=\{u \in L^2:\ \mathcal E_K[u]<\infty\} = \Lambda_{2,2}^{\beta/2}$, where $\beta = \log \lambda / \log r_0$.
\end{theorem}

\medskip

We set $\Vert u \Vert_{\mathcal E_K} = \Vert u \Vert_{L^2} + (\mathcal E_K[u])^{1/2}$. Then $(\mathcal D_K, \Vert \cdot \Vert_{\mathcal E_K})$ is a Banach space which is equivalent to $(\Lambda_{2,2}^{\beta/2}, \Vert\cdot\Vert_{\Lambda_{2,2}^{\beta/2}})$. By the polarization identity,  the energy form defines a symmetric bilinear form.  We  first consider a situation where $(\mathcal D_K, \Vert \cdot \Vert_{\mathcal E_K})$ consists of  H\"older continuous functions.

\medskip

Let ${\mathcal C}^\delta$ be the space of H\"older continuous functions on $K$ of order $\delta$ and let
$\Vert u\Vert_{ {\mathcal C}^\delta}:= \Vert u\Vert_\infty + {\rm sup}_{\xi,\eta \in K} \frac{|u(\xi)-u(\eta)|}{\rho(\xi,\eta)^\delta}$.
We define the {\it upper dimension} $\bar d_\mu = \overline{\dim}\,\mu$ of $\mu$ on $(K,\rho,\mu)$ by
\begin{equation} \label{eq5.1}
\bar d_\mu =  \inf\{\alpha :  \exists \ c>0 \text{ s.t. } \mu(B(\xi,r)) \geq c r^\alpha \ \forall \ \xi \in K \hbox{ and } 0<r\leq 1\}.
\end{equation}
(Recall that $\diam (K) =1$ by convention.)
Using \eqref{eq5.1}, it is direct to extend the embedding theorem in \cite{GHL1} to the following.

\medskip

\begin{proposition} \label{th5.1}
Let $(K,\rho,\mu)$ be a compact space of homogeneous type. If $\beta>\alpha>\bar d_\mu$, then
\begin{equation*}
\Vert u \Vert_{\mathcal C^{(\beta-\alpha)/2}} \leq C \Vert u \Vert_{\Lambda_{2,2}^{\beta/2}}, \qquad \forall \ u \in L^2.
\end{equation*}
That is, $\Lambda_{2,2}^{\beta/2} \hookrightarrow \mathcal C^{(\beta-\alpha)/2}$ is an embedding.
\end{proposition}

\bigskip
%

Assume  the $\lambda$-NRW on the augmented tree of $(K, \rho, \mu)$  has a return ratio $\lambda \in (0,r_0^{\bar d_\mu})$, that is, $\beta = \log \lambda / \log r_0 > \bar d_\mu$. Then for any $\alpha \in (\bar d_\mu, \beta)$ fixed, we have $\mathcal D_K = \Lambda_{2,2}^\beta \subset {\mathcal C}^{(\beta-\alpha)/2} \subset  C(K)$ by Proposition \ref{th5.1}, and $c(x,x^-) \gtrsim (r_0^\alpha/\lambda)^{|x|}$ for all $x \in X$ by \eqref{eq4.6} and \eqref{eq5.1}. It follows that
\begin{equation} \label{eq5.3}
|f(x)-f(y)| \leq C (\lambda/r_0^\alpha)^{|x|/2}\mathcal E_X[f]^{1/2}, \qquad \forall \ f \in \mathcal D_X \hbox{ and } (x,y) \in E,
\end{equation}
and therefore we can define a {\it trace map} ${\rm Tr}  : \mathcal D_X \to \mathcal D_K $ by
\begin{equation} \label{eq5.4}
({\rm Tr} f)(\xi) = \lim_{n \to \infty} f(x_n), \qquad \xi \in K \hbox{ and }\pi(\vartheta,\xi) = [\vartheta,x_1,\cdots].
\end{equation}
Clearly the above limit is independent of the choice of $\pi(\vartheta,\xi)$ and is uniform on $K$.
Next recall that the Poisson integral $H$ in \eqref{eq2.6} maps $\mathcal D_K$ into $\mathcal{HD}_X$, the class of harmonic functions in $\mathcal D_X$. We further impose a complete norm $\Vert \cdot \Vert_{\mathcal E_X}$ on $\mathcal D_X$ by
\begin{equation*}
\Vert f \Vert_{\mathcal E_X}^2 = {\sum}_{x \in X}|f(x)|^2 w^{|x|}+\mathcal E_X[f],
\end{equation*}
 where $0 < w < r_0^{\bar d_\mu}$ is a constant  (see \cite[Corollary 3.6]{KL}). Applying Corollary \ref{th4.5'} and \eqref{eq5.3}, we obtain the following trace theorem \cite[Section 3]{KL}.

\medskip

\begin{theorem} \label{th5.3} \hspace{-2mm}
Let $(K,\rho,\mu)$ be a compact space of homogeneous type, and let $\{Z_n\}_{n=0}^\infty$ be a $\lambda$-NRW with ratio $\lambda \in (0,r_0^{\bar d_\mu})$ on the augmented tree $(X,E,\Phi)$. Then ${\rm Tr}(\mathcal D_X) = \mathcal D_K$, and ${\rm Tr}|_{\mathcal{HD}_X} = H^{-1}$. Moreover, ${\rm Tr}: (\mathcal{HD}_X, \Vert\cdot\Vert_{\mathcal E_X}) \to (\mathcal D_K, \Vert\cdot\Vert_{\mathcal E_K})$ is a Banach space isomorphism.
\end{theorem}

We are interested in conditions under which $(\mathcal E_K,\mathcal D_K)$ is a (non-local) Dirichlet form. It is easy to check the properties in
Definition \ref{de5.1} except the density property in $L^2$ and in $C(K)$.  To this end, we introduce two {\it critical exponents} of the family $\{\Lambda_{2,2}^{\beta/2}\}_{\beta>0}$ by
\begin{align*}
\beta^* &:= \sup \{\beta > 0: \Lambda_{2,2}^{\beta/2} \cap  C(K) \hbox{ is dense in } C(K) \hbox{ with } \Vert\cdot\Vert_\infty\}, \quad \hbox{and} \\
\beta^\sharp &:= \sup\{\beta > 0: \dim(\Lambda_{2,2}^{\beta/2} \cap  C(K)) > 1\}.
\end{align*}
As $C(K)$ is dense in $L^2(K, \rho, \mu)$, the Proposition \ref{th5.3} implies that $(\mathcal E_K, \mathcal D_K)$ is a non-local regular Dirichlet form if $\bar d_\mu < \beta < \beta^*$.
In general  when $\rho$ is a metric, we know that \ $2 \leq \beta^* (\leq \beta^\sharp \leq \infty)$ (see \cite{St} for $\alpha$-set $K$).  Indeed, if we let ${\mathcal G}:= \{\rho_\xi\}_{\xi \in K}$ be the class of distance functions on $K$ defined by $\rho_\xi(\eta) := \rho(\xi,\eta)$. Then it is straightforward that ${\mathcal G}\subset \Lambda_{2,2}^{\beta/2} \cap  C(K)$ for $\beta \in (0,2)$, and ${\mathcal G}$ {\it separates points} in $K$.  By the Stone-Weierstrass theorem, $\Lambda_{2,2}^{\beta/2} \cap  C(K)$ is dense in $C(K)$ when $\beta \in (0,2)$, and hence $2 \leq \beta^* $. (For $\beta \in (0,2)$, we  cannot prove whether $\Lambda_{2,2}^{\beta/2} \cap C(K)$ is dense in $\Lambda_{2,2}^{\beta/2}$ with $\Vert\cdot\Vert_{\mathcal E_K}$.) Furthermore,  by using the same argument as in \cite{GHL1}, we can show that $\beta^*$ has an upper bound, $\beta^* \leq \beta^\sharp \leq \bar d_\mu + 1$ provided that $(K,\rho)$ satisfies a {\it chain condition} in \cite{GHL1}.

\medskip

The values of $\beta^*$ and $\beta^\sharp$ are known and are equal for any standard cases: for  any  classical domain in $\mathbb R^d$ with Lebesgue measure, $\beta^*=\beta^\sharp=2$; for the $d$-dimensional Sierpi{\'n}ski gasket with $\alpha$-Hausdorff measure (here $\alpha = \log (d+1)/2$ is the Hausdorff dimension), $\beta^*=\beta^\sharp=\log (d+3)/\log 2$ (see \cite{Jo});  for Cantor-type set, $\beta^*=\beta^\sharp=\infty$. There are also examples that the two exponents $\beta^*$ and $\beta^\sharp$ are different (see \cite{GuL,KL}).

\bigskip

Corresponding to the critical exponents of the domain ${\mathcal D}_K$ of the induced energy $\mathcal E_K$, there are critical values of the return ratio $\lambda$ of the random walk on the augmented tree $(X,E)$. Such values can be analyzed, and explicitly calculated in some special cases through the associated  electrical network on $(X,E)$. For this, we recall  the  effective resistances defined as in \eqref{eq2.4}; we will fix a set of conductance that defines a $\lambda$-NRW (it can be replaced by any other one in \eqref{eq4.6}) by letting
\begin{equation} \label{eq5.6}
c(x,x^-) = \lambda^{-|x|}\mu(\Phi(x)), \ \ x \in X \setminus \{\vartheta\}, \quad  c(x,y) = \frac{c(x,x^-)c(y,y^-)}{c(x,x^-)+c(y,y^-)}, \ \  (x,y) \in E_h.
\end{equation}
 Let $\{\kappa_n\}_{n=0}^\infty$ be a sequence of mappings from $K$ to $X$ such that for each $\xi \in K$, the sequence $[\kappa_0(\xi), \cdots , \kappa_n(\xi) , \cdots]$ is a geodesic ray $\pi(\vartheta,\xi)$; we call such $\{\kappa_n\}_{n=0}^\infty$ a {\it $\kappa$-sequence}. For $n \geq 1$, let $X_n:=\bigcup_{k=0}^n \mathcal J_k$, and define the {\it level-$n$ resistance} (with respect to $\kappa_n$) between two nonempty closed subsets $A, B \subset K$ by
\begin{equation} \label{eq5.7}
R_n^{(\lambda)}(A,B):= R_{X_n}(\kappa_n(A),\kappa_n(B)),
\end{equation}
where the conductance on $X_n \times X_n$ is given by \eqref{eq5.6}. With the assumption $\lambda \in (0,r_0^{\bar d_\mu})$, it has been proved that the limit
\begin{equation} \label{eq5.8}
R^{(\lambda)}(A,B):= \lim_{n \to \infty} R_n^{(\lambda)}(A,B)
\end{equation}
exists, and is independent of the choice of the $\kappa$-sequence $\{\kappa_n\}$ \cite[Theorem 4.2]{KL}; we call $R^{(\lambda)}(A,B)$ the {\it limit resistance} between $A$ and $B$. For $\xi,\eta \in K$, we  write $R^{(\lambda)}(\{\xi\},\{\eta\}) = R^{(\lambda)}(\xi,\eta)$.

\medskip

\begin{theorem} \label{th5.4} \hspace{-2mm} {\rm \cite[Theorem 5.1 and  Corollary 5.2]{KL}} Let $(K,\rho,\mu)$ be a compact space of homogeneous type, and let $\{Z_n\}_{n=0}^\infty$ be a $\lambda$-NRW with ratio $\lambda \in (0,r_0^{\bar d_\mu})$ on the augmented tree $(X,E,\Phi)$. Suppose two nonempty closed subsets $A, B \subset K$ satisfy $R^{(\lambda)}(A,B)>0$. Then
\begin{equation} \label{eq5.9}
R^{(\lambda)}(A, B)^{-1} = \inf \{\mathcal E_K[u]: u \in \mathcal D_K \hbox{ with } u = 1 \hbox{ on } A, \hbox{ and } u=0 \hbox{ on } B\},
\end{equation}
and there exists a unique energy minimizer $u_0$, i.e., $u_0|_A = 1$, $u_0|_B = 0$ and $\mathcal E_K[u_0]=R^{(\lambda)}(A,B)^{-1}$.
In particular, for $\xi,\eta \in K$, $R^{(\lambda)}(\xi,\eta)>0$ if and only if there exists $u \in \mathcal D_K$ such that $u(\xi) \neq u(\eta)$; in this case, we have
$$
R^{(\lambda)}(\xi,\eta) = \sup \left\{\frac{|u(\xi)-u(\eta)|^2}{\mathcal E_K[u]}: u \in \mathcal D_K, \ \mathcal E_K[u]>0 \right\}.
$$
\end{theorem}

\bigskip
We introduce two critical values for the limit resistances by
\begin{align*}
\lambda^* &:= \inf\{\lambda \in (0,r_0^{\bar d_\mu}): R^{(\lambda)}(\xi,\eta)>0, \ \forall \ \xi \neq \eta \in K\}, \quad \hbox{and} \\
\lambda^\sharp &:= \sup\{\lambda \in (0,r_0^{\bar d_\mu}): R^{(\lambda)}(\xi,\eta)=0, \ \forall \ \xi, \eta \in K\}.
\end{align*}
Then $\lambda^\sharp \leq \lambda^*$, and the above theorem together with the Stone-Weierstrass theorem yield

\medskip

\begin{corollary} \label{th5.5}
With the same assumption as in Theorem \ref{th5.4}, if $\lambda^* (\lambda^\sharp )\in (0,r_0^{\bar d_\mu})$, we have
$$
\beta^* = \log \lambda^* / \log r_0, \quad ( \,  \beta^\sharp = \log \lambda^\sharp / \log r_0, \hbox { respectively}),
$$
and if $\lambda^* = 0$, then $\beta^*= \infty$.
\end{corollary}

\medskip
%

In general when the cardinality of $K$ is infinity, it is difficult to determine the critical exponents $\beta^*$ and $\beta^\sharp$ by applying Corollary \ref{th5.5} directly, since it is not feasible to test the positivity of limit resistances between every pair of points in $K$. However, when $K$ is a self-similar set, such an infinite exhaustive testing can be reduced to a finite set.

\medskip

For a self-similar set $K$, we use the same notations as in Example \ref{ex3.3}. Let $i^\infty \in K$ denote the unique fixed point of the contractive similitude $S_i$ in the IFS $\{S_i\}_{i \in \Sigma}$, i.e., $\{i^\infty\} = \bigcap_{n=1}^\infty {S_{i^n}(K)}$.

\medskip

\begin{theorem} \label{th5.6} \hspace{-2mm} {\rm \cite[Theorem 5.4]{KL}}
Suppose the IFS $\{S_i\}_{i \in \Sigma}$ satisfies the OSC. Let $K$ be the self-similar set equipped with a doubling self-similar measure $\mu$, and let $(X,E,\Phi)$ be the augmented tree defined in Example \ref{ex3.3} and \eqref{eq3.3}. Then
\begin{equation*}
\lambda^\sharp = \sup\{\lambda \in (0,r_0^{\bar d_\mu}): R^{(\lambda)}(i^\infty, j^\infty)=0, \ \forall \ i,j \in \Sigma\}.
\end{equation*}
\end{theorem}

\bigskip

We say that a self-similar set $K$ is {\it post critically finite} (p.c.f.) \cite{Ki1} if it has a finite {\it post critical set} $\mathcal P$ defined by
$$
\mathcal P = {\bigcup}_{n \geq 1} \sigma^n \Big(\pi^{-1}\big({\bigcup}_{i,j \in \Sigma, i \neq j} (S_i(K) \cap S_j(K))\big)\Big),
$$
where $\pi$ is the natural projection from the symbolic space $\Sigma^\infty$ to $K$, and $\sigma$ is the shift operator on $\Sigma^\infty$. We introduce a geometric condition on the p.c.f.~set:

\vspace{2mm}

\noindent {\rm ($\ast$)} \ {\it there exist $\delta, c>0$ such that for any $i,j \in \Sigma$ and $\zeta \in S_i(K) \cap S_j(K)$,}
\begin{equation*}
|\xi-\eta| \geq c(|\xi-\zeta|+|\zeta-\eta|), \qquad \forall\ \xi \in S_i(K) \cap B(\zeta,\delta) \hbox{ and } \eta \in S_j(K) \cap B(\zeta,\delta).
\end{equation*}
It is easy to check that the most familiar p.c.f.~sets including nested fractals satisfy this condition. Let $V_0 = \pi(\mathcal P)$ be the ``boundary" of a p.c.f.~set $K$.

\begin{theorem} \label{th5.7} \hspace{-2mm} {\rm \cite[Theorem 5.9]{KL}}
With the same assumption as in Theorem \ref{th5.6}, assume further $K$ is p.c.f.~and satisfies ($\ast$). Then
\begin{equation*}
\lambda^* = \inf\{\lambda \in (0,r_0^{\bar d_\mu}): R^{(\lambda)}(\xi,\eta)>0, \ \forall \ \xi \neq \eta \in V_0\}.
\end{equation*}
\end{theorem}

\medskip

In \cite{KL} we provided a procedure and examples to implement the above two theorems; the technique  is from the  network theory including  the {\it series} and {\it parallel laws},  the {\it monotonicity law} ({\it cutting} and {\it shorting}) as well as the  {\it $\Delta$-Y transfrom} \cite{DS,LyP}. We remark that on the same underlying set, there can be different augmented trees which admit different random walks and induces different Dirichlet forms. For example for $K=[0,1]$, it is generated by $S_1(\xi) = \frac 12 \xi$ and $S_2(\xi) = \frac 12 (\xi+1)$ on $\mathbb R$. The self-similar measure $\mu$ with doubling property \eqref{eq3.1} can only be the Lebesgue measure; in this case, $\beta^*=\beta^\sharp = 2$. On the other hand, by using the same method as in \cite{KL}, we see in the following example that for $[0,1]$ generated by another iterated function system, the results can be vastly different.

\medskip

\begin{example} \label{ex5.8} {\bf Rotated unit interval} {\rm Let $S_1(\xi) = \frac 12 \xi$ and $S_2(\xi) = 1-\frac 12 \xi$ on $\mathbb R$. Then the self-similar set $K$ is the interval $[0,1]$, and is p.c.f.~with $\mathcal P=\{1^\infty,21^\infty\}$. Let $\mu$ be the self-similar measure generated by probability weights $p, q$. Then $(K,  \mu)$ is doubling with respect to the Euclidean metric on $\mathbb R$ (see \cite{Y}), and its upper dimension is $\bar d_\mu = -\log (\min\{p,q\}) / \log 2$. The $\lambda^* = \lambda ^\sharp = pq$,  and the critical exponents $\beta^* = \beta^\sharp = -\log(pq)/\log 2$.}
\end{example}

\medskip

\bigskip
To conclude, we remark that the induced Dirichlet forms ${\mathcal E}_K^{(\beta)}$ in our study are non-local. For local regular Dirchlet forms, which give the Laplacians, are  more difficult to obtain. So far in the analysis on fractals, the most studied Laplacians are on the p.c.f. sets and the Sierpi{\'n}ski carpet. It is still an open question for the existence of local Dirichlet form/Laplacian on  self-similar sets; in particular, we do not know if there is a non-trivial fractal set that do not support a Laplacian. On the other hand, it is observed that in all the known cases, the local regular Dirichlet forms have Besov spaces  $\Lambda^{\beta^*}_{2, \infty}$ as domains, and $\Lambda^{\beta^*}_{2, 2}$ only contain constant functions. The critical exponent conceals  lot of information, and it is still far from clear. One of the challenge questions is  to understand the  ``transition" of ${\mathcal E}_K^{(\beta)}$  as $\beta \nearrow  \beta^*$ (see \cite {BBM, GY, Ya}), and in addition, in terms of the random walks in our consideration.

\bigskip
\small{
}

\bigskip
\bigskip

\noindent Shi-Lei Kong, Fakult{\"a}t f{\"u}r Mathematik, Universit{\"a}t Bielefeld, Postfach 100131, 33501 Bielefeld, Germany. \\
skong@math.uni-bielefeld.de

\bigskip

\noindent Ka-Sing Lau, Department of Mathematics, The Chinese University of Hong Kong, Hong Kong.\\
\& School of Mathematics and Statistics, Central China Normal University, Wuhan, 430079, China. \\
\& Department of Mathematics, University of Pittsburgh, Pittsburgh, PA 15260, USA. \\
kslau@math.cuhk.edu.hk

\bigskip

\noindent Ting-Kam Leonard Wong, Department of Mathematics, University of Southern California, Los Angeles, CA 90089 USA.\\
tkleonardwong@gmail.com

\end{document}